\documentclass[11pt, a4paper]{amsart}

\usepackage{amsmath}
\usepackage{latexsym}
\usepackage[all]{xy}
\usepackage{color}
\newtheorem{teorema}{Theorem}[section]
\newtheorem{definicion}[teorema]{Definition}
\include{amslatex}
\newtheorem{proposicion}[teorema]{Proposition}

\newtheorem{ejemplo}[teorema]{Example}

\numberwithin{equation}{section}

\begin{document}
\title[Average in Finsler spaces of Lorentzian structure]{Average of geometric structures in Finsler spaces with Lorentzian signature}

\maketitle
\begin{center}
\author{Ricardo Gallego Torrom\'e\footnote{email: rigato39@gmail.com}}
\end{center}
\begin{center}
\address{Department of Mathematics\\
Faculty of Mathematics, Natural Sciences and Information Technologies\\
University of Primorska, Koper, Slovenia}
\end{center}

\begin{abstract}
Given the class of Finsler spaces with Lorentzian signature $(M,L)$ on a manifold $M$ endowed with a timelike vector field $\mathcal{X}$ satisfying $g_{(x,y)}(\mathcal{X},\mathcal{X})<0$ at any point $(x,y)$ of the slit tangent bundle, a pseudo-Riemannian metric defined on $M$  of signature $n-1$ is associated to the fundamental tensor $g$. Furthermore, an affine, torsion free connection is associated to the Chern connection determined by $L$. The definition of the average connection does not make use of $\mathcal{X}$. Therefore, there is no direct relation between these two averaged objects.
\end{abstract}
\bigskip
{\bf MSC Class}: 53DC60; 53B40; 83D05.

\section{Introduction}

Positive definite Finsler spaces are natural generalizations of Riemannian spaces such that for the associated norm function, the Riemannian assumption of quadratic dependence on the coordinate velocities is dropped \cite{Chern}. The theory of positive definite Finsler spaces has been developed in the form of sound geometric theories that encompasses the Riemannian case if the above mentioned quadratic restriction in the norm is imposed \cite{Cartan Exposes II, Berwald, Chern1948}, geometric frameworks where fundamental results from Riemannian geometry are conveniently generalized to the Finslerian case \cite{BCS}.

For positive definite signatures, the theory of average Finsler structures provides a map from the Finsler category to the Riemannian or affine categories, depending on the object considered, by conveniently {\it tracing out} the $y$-dependence of relevant Finslerian geometric object. Apart from the general form of the theory discussed in \cite{Ricardo05}, other theories of averaging Finsler metrics have appeared in the literature. For instance, in the context of Szabo's metrizability theorem \cite{Szabo1980}, a theory of averaging the fundamental tensor of a Finsler structure was discussed in the work of  Vincze in \cite{Csava2005} and, in a different formulation, the averaging the Finsler metric is fundamental for the methods applied by  Matveev and Troyanov  in \cite{Matveev Troyanov 2012}.

The sister theory of Finsler spaces with Lorenztian structures appears naturally in the form of relativistic Finslerian generalizations of Lorentzian spacetimes \cite{Ricardo 2017b}. In contrast with the positive case, for Lorentzian signature there is no currently an unique theoretical paradigm where all the questions and specific models found in the literature can be conveniently formulated. Indeed, several theories of Finsler structures with Lorentzian signature have been developed in the literature (see \cite{Pfeifer 2019} for a review of different recent approaches), the main problems related with the appearance of singularities for the fundamental tensor and the impossibility to treat within the same formalism relevant particular cases.

In the  context described above, one still can consider the problem of generalizing the average operation method to Finsler spaces with Lorentzian signature. Due to the naturalness of the averaging methods in the positive definite case and the interpretation of coarse grained operation in the Lorentzian signature case, this problem is of interest from the point of view of differential geometry and also for its applications in extensions of general relativity. However, the direct adaptation of the procedure as in positive Finslerian geometry leads to singularities of the average metric where the original Finsler metric was regular, a problem that spoils the construction of significant connections from the average metric.

 In this work we consider how to define a meaningful averaging operation for the fundamental tensor and for the Chern connection determined by a Finsler structure with Lorentzian signature. This question is investigated in the framework of J. Beem theory \cite{Beem1} for Finsler spaces with Lorentzian signature. For this purpose, we have restricted the attention to metric structures allowing the existence of a vector field $\mathcal{X}$ timelike respect to any direction of the tangent space (see the relation \eqref{condition on X}). Such condition is stronger than the condition of time orientation in the Finslerian sense \cite{GallegoPiccioneVitorio:2012}, but can be motivated physically as we will discuss below. On the other hand, the definition of the average Chern connection for Finsler spaces of Lorentzian signature is more direct, but it faces the problem of the non-compactness of the unit hyperboloid. Besides this complication, the formulation of the average connection for Finsler metrics with Lorentzian signature is a direct generalization of the positive definite case.

The structure of this work is the following. In section 2 we introduce the necessary notions of Finsler geometry, including the averaging method as discussed in reference \cite{Ricardo05}. The fundamental concepts of Finsler spaces with Lorentzian signature that we will need are also introduced, according to \cite{Beem1, GallegoPiccioneVitorio:2012}. In section 3, we introduce our definition of average Finsler metric for structures of Lorentzian signature. The method associates to the initial  Finsler structure with Lorentzian signature a positive definite Lagrange metric \cite{MironHrimiucShimadaSabau:2002}. The averaging operation is applied to such positive Lagrange metric, obtaining a Riemannian metric on $M$. After that, a metric of signature $n-1$ defined over $M$ by the method discussed in conventional Lorentzian geometry (see for instance \cite{Hawking Ellis 1973}, section 2.6) is obtained. It is only in the first step of this process that the condition \eqref{condition on X} on $\mathcal{X}$ is necessary and used.
Section 4 deals with the average of the Chern connection. Differently than in the positive case, in Lorentzian signature the non-compactness of the unit hyperboloid introduces an extra-difficulty. Two ways to overcome this difficulty are briefly discussed.

\section{Geometric framework}
\subsection{Geometric framework for average of Finsler metrics and connections}
We introduce the average operation of Finsler metrics and connections following \cite{Ricardo05}. Let $M$ be a differentiable manifold of dimension $n=\,dim( M)$ and $TM$ the tangent bundle of $M$. Local
coordinates $(U,x)$ on $M$ induce local natural coordinates $(TU, x^i, y^i)$ on $TM$. The slit tangent bundle is $N=TM\setminus\{0\}$,
where $0$ is the zero section of $TM$; $N_x\subset\,N$ is the fiber over $x\in \,M$. $\pi^*TM\to N$ will denote the pull-back bundle.
\begin{definicion}
\label{Finslerstructute}
A Finsler space on the manifold $M$  is a pair $(M,F)$ where $F$ is a
non-negative, real function  $F:M\to [0,+\infty [$
such that
\begin{itemize}
\item It is smooth in the slit tangent bundle $N$,

 \item Positive homogeneity holds: $F(x,{\lambda}\,y)=\,\lambda \,F(x,y)$ for every $\lambda \in\,[0,+\infty [$,

 \item Strong convexity holds:
the Hessian matrix
 \begin{align}
g_{ij}(x,y): =\,\frac{1}{2}\,\frac{{\partial}^2 F^2
(x,y)}{{\partial}y^i {\partial}y^j },\quad i,j=\,1,...,n
\label{fundamentaltensorcoefficients}
\end{align}
is positive definite on $N$.
\end{itemize}
\label{definicionfinslerstructure}
\end{definicion}
The $n$-form $d^n y$ is defined in local coordinates by the expression
\begin{align}
d^n y=\sqrt{|\det\, g(x,y)|}\,\, \delta y^1 \wedge \cdot \cdot \cdot \wedge \delta y^n,
\label{volumeformontangentspace}
\end{align}
where $\det\, g(x,y)$ is the determinant of the fundamental tensor
\begin{align}
\left(g\right)_{ij}(x,y):=g_{ij}(x,y),\quad i,j=1,...,n.
\label{fundamental tensor}
 \end{align}
 and $\{\delta y^1,...,\delta y^n\}$ is a local dual frame to the holonomic frame $\{\frac{\partial}{\partial y^1},...,\frac{\partial}{\partial y^n}\}$. $\{\delta y^1,...,\delta y^n\}$  is constructed using the Finsler norm $F$. The details can be found in \cite{BCS} for the positive case and in \cite{GallegoPiccioneVitorio:2012} for the case of signature $n-1$.
 Note that in \cite{BCS}, there are certain factors in terms of powers of $F$  in the definition of the above quantities respect to our definitions. Such factors are put on hand to meet  certain homogeneous properties for the local basis and geometric quantizes. Such factors do not appear in the Lorentzian case, since in this case the condition $F=0$ can lead to un-necessary apparent singularities in the local expression of geometric objects \cite{GallegoPiccioneVitorio:2012}.
For a positive definite Finsler space, the indicatrix  over $x$ is the compact submanifold $I_x :=\{y\in T_xM\,s.t.\,F(x,y)=1\}\hookrightarrow \,T_xM$.

Since $d^n y$ is invariant by local diffeomorphisms on $TM$, it defines a section of $\Lambda^n N$.
 For each embedding $i_x:I_x \hookrightarrow N_x$ one consider the volume form on $I_x$ given by
\begin{align}
dvol_x:=i^*_x (d^n y\cdot\l ),
\label{volumeformontheindicatrix}
\end{align}
where $\l =\,{y^i}\,\frac{\partial}{\partial y^i}$ and $d^n y\cdot\l$ is the corresponding contraction; $i^*_x$ is the pull-back on $T^*N_x$.
Then the volume function $vol(I_x)$ is defined by the expression
\begin{align}
vol:M\to \mathbb{R}^+,\quad x\mapsto vol(I_x)=\int _{I_x} \, |\psi|^2(x,y)\,dvol_x,
\label{volumenoftheindicatrix}
\end{align}
where the weight factor $|\psi|^2:TM\setminus\{0\}\to \mathbb{R}^+$ is an homogenous of degree zero in $y$, positive, smooth function on the tangent bundle. Let $\pi_{\mathcal{I}}:\mathcal{I}\to M$ be the indicatrix bundle over $M$, where each fiber $\pi^{-1}_{\mathcal{I}}(x)=\,I_x$.
The average of a function $f\in \mathcal{F}(\mathcal{I})$ is defined by the expression
\begin{align*}
\langle f\rangle_{\psi} (x):=\frac{1}{vol(I_x)} \int _{I_x }\,|\psi|^2(x,y)\,f(x,y)\,dvol_x.
\end{align*}

Given the Finsler space $(M,F)$, let us consider the matrix with components
\begin{align}
h_{ij}(\psi,x):=\langle g_{ij}(x,y)\rangle_{\psi},\quad i,j=1,...,n,
\label{averagedmetric}
\end{align}
for each point $x\in M$.
\begin{proposicion}
Let $(M,F)$ be a Finsler space. Then $\{h_{ij}(\psi, x)\}^{n}_{i,j=1}$
are the components of a Riemannian metric
\begin{align}
h_\psi(x)=h_{ij}(\psi,x)\, dx^i\otimes dx^j,\quad i,j=1,...,n.
\label{averagedmetric2}
\end{align}
\label{propositionaveragemetric}
\end{proposicion}
The proof of this result can be found in \cite{Ricardo05}. This is the fundamental result on which the notion average Finsler metric is built on and that we will generalize to certain metric structures with signature $n-1$, although the methodology of the positive case cannot be transferred directly to the Lorentzian signature case.

 An analogous construction can be applied to the connections determined by a Finsler structure, for instance, the Chern's connection, the Cartan's connection and other geometric objects. The details of these constructions can be found in \cite{Ricardo05}. We collect here the result that will be considered for generalization to the Lorentzian case. Let $\pi^*TM\to M$ be the pull-back vector bundle of $\pi_N:N\to M$ and $u\in \,N_x$. Also, let us remind the reader that it is possible to define a non-linear connection on $N$ from the Finsler metric. Given such non-linear connection, $h_x :T_xM\to T_uN$ is the horizontal lift map \cite{Ricardo05}. Then we have the following
\begin{teorema}
 Let $\nabla$ be a linear connection of the vector bundle $\pi^*TM\to N$. Then there
is defined an affine connection on $M$ determined by the
covariant derivative of each section $Y\in\,\Gamma\, TM$ along each  direction $X\in T_x M$,
\begin{align}
\langle\nabla\rangle_X  Y:= \langle \pi _2|_u
{\nabla}_{h_u({X})} \pi^* _v Y\,\rangle ,\,v\in TU_x\setminus \{0\},
\label{averagedconnectionformula}
\end{align}
for each $X\in T_x M$ and $Y\in\,\Gamma\, TM$, where ${U}_x$ is an open neighborhood of $x\in M$.
\label{averagedconnection}
\end{teorema}
This theorem defines the averaged Chern connection.
In adequate local frames, the connection coefficients of the averaged connection are the average of the connection coefficients of the Chern connection.
\subsection{Geometric framework for Finsler spaces with Lorentzian signature}
Following  J. Beem formalism \cite{Beem1} but using the notation from \cite{GallegoPiccioneVitorio:2012},
we introduce the basic notation and fundamental notions of Finsler spaces with Lorentzian signature theory. As we mentioned before, we formulate the theory for $n$-dimensional manifolds. Thus we start with the following
\begin{definicion}
A Finsler spaces with Lorentzian signature is a pair $(M,L)$ where
\begin{enumerate}
\item $M$ is an $n$-dimensional real, second countable, Hausdorff $C^{\infty}$-manifold.
\item $L:N\longrightarrow R$ is a real smooth function such that
\begin{enumerate}
\item $L(x,\cdot)$ is positive homogeneous of degree two in the variable $y$,
\begin{align}
L(x,ky)=\,k^2\,L(x,y),\quad \forall\, k\in ]0,+\infty[,
\label{homogeneouscondition}
\end{align}
\item The {\it vertical Hessian}
\begin{align}
g_{ij}(x,y)=\,\frac{\partial^2\,L(x,y)}{\partial y^i\,\partial y^j}
\label{nondegeracy-signature}
\end{align}
is non-degenerate and with signature $n-1$ for all $(x,y)\in\, N$.
\end{enumerate}
\end{enumerate}
\label{Finslerspacetime}
\end{definicion}

Direct consequences of this definition and Euler's theorem for positive homogeneous functions are the following relations,
\begin{align}
\frac{\partial L(x,y)}{\partial y^k}\,y^k=\,2\,L(x,y),\quad \frac{\partial L(x,y)}{\partial y^i}=\,g_{ij}(x,y)y^j,
\quad L(x,y)=\frac{1}{2}\,g_{ij}(x,y)y^iy^j.
\end{align}
Furthermore, the function $L$ determines an homogeneous of  degree one function $F$ by the relation $F:=L^{1/2}$. One main difference with the positive definite Finsler case is that the function $F$ is not real valued in the whole $TM$.

Given a Finsler spaces with Lorentzian signature $(M,L)$,
a vector field $X\in\,\Gamma TM$ is timelike if $L(x, X(x))<0$ for every point $x\in\,M$ and a curve $\lambda:I\longrightarrow M$ is timelike
if the tangent vector field is timelike in the sense that $L(\lambda(s),\dot{\lambda}(s))<0$. A vector field $X\in \,
 \Gamma TN$ is lightlike if $L(x,X(x))=0$ for every point $x\in\,M$ and a curve is lightlike if its tangent
  vector field is lightlike. Similar notions are for spacelike.
A curve is causal if either it is timelike and has constant
  speed $g_{\dot{\lambda}}(\dot{\lambda},\dot{\lambda}):=L(\lambda,\dot{\lambda})=g_{ij}(\dot{\lambda},T)\dot{\lambda}^i\dot{\lambda}^j$ or if it is lightlike.
There is a class sub-manifolds, one on each tangent space, that play a significant role in the theory. The {\it unit hyperboloid} at $x\in M$ is defined by
\begin{align*}
\Sigma_x:=\{y\in\,T_x M\,st.\,L(x,y)=-1\}.
 \end{align*}
$\Sigma_x$ is a co-dimension $1$ sub-manifold of $N_x$. 

The Cartan tensor is defined here as a tensor field of components
\begin{align}
C_{ijk}:=\,\frac{1}{2} \frac{\partial g_{ij}}{\partial y^{k}}, \quad i,j,k=1,...,n,
\label{Cartantensor coefficients}
\end{align}
slightly different than in the way it is formulated in positive definite Finsler geometry \cite{BCS, Cartan Exposes II}. As mentioned before, we do not use factors of $F$ to make the Cartan tensor homogeneous of degree zero. However, the fundamental tensor $g$ is still of 0-homogeneity in $y$-variables. Therefore, by Euler's theorem
\begin{align}
C_{(x,y)}(y,\cdot,\cdot)=\,\,\frac{1}{2} y^k\,\frac{\partial g_{ij}}{\partial y^{k}}=0.
\label{homogeneityofcartan tensor}
\end{align}

Chern connection for Finsler spaces with Lorentzian signature is defined in similar terms as in the case of positive definite Finsler spaces. We introduce here the connection following the index-free formulation that can be found in \cite{GallegoPiccioneVitorio:2012},
\begin{proposicion}
Let $h(X)$ and $v(X)$ be the horizontal and vertical lifts of $X\in\,\Gamma TM$ to $TN$, and $\pi^*g$ the pull back-metric.
The Chern connection is characterized by
\begin{enumerate}
\item  The almost $g$-compatibility metric condition is equivalent to
\begin{align}
{\nabla}_{v({\hat{X}})}\pi^*g=2C(\hat{X},\cdot,\cdot),\quad
{\nabla}_{h({X})} \pi^*g=0,\quad \hat{X}\in\,\Gamma TN.
\label{covariantalmostmetriccondition}
\end{align}

\item The torsion-free condition:
\begin{enumerate}
\item Null vertical covariant derivative of sections of ${\pi}^*{TM}$:
\begin{align}
{\nabla}_{{V({X})}} {\pi}^* Y=0,
\label{covariantensorfreecondition1}
\end{align}
for any vertical component $V(X)$ of $X$.
\item  Let us consider $X,Y\in {TM}$ and their horizontal
lifts $h(X)$ and $h(Y)$. Then
\begin{align}
\nabla_{h(X)} {\pi}^* Y-{\nabla}_{h(Y)}{\pi}^*X-{\pi}^* ([X,Y])=0.
\label{covariantensorfreecondition2}
\end{align}
\end{enumerate}
\end{enumerate}
\end{proposicion}
The connection coefficients $\Gamma^i_{jk}(x,y) $ of the Chern connection are constructed  in terms of the fundamental tensor components and the Cartan tensor components. The detail of how they are constructed in the setting of Finsler spacetimes can be found in detail in \cite{GallegoPiccioneVitorio:2012}.

\section{Definition of the average of the metric tensor for Finsler spaces with Lorentzian signature}
The direct extension to Lorentzian signature of the average metric tensor for positive definite Finsler metrics by means of  \eqref{averagedmetric2}  is not a good strategy, because the average of the components $\langle g_{ij}(x,y)\rangle_\psi$ does not provide a well-defined, non-singular metric tensor.
Instead, the procedure that we follow for the generalization requires of three steps and is the following.
\\
{\bf 1}. Let us assume that $M$ is endowed with a vector field $\mathcal{X}\in\,\Gamma \,TM$ such that
\begin{align}
g_{(x,y)}(\mathcal{X},\mathcal{X})<0,\,\quad \forall \,(x,y)\in\,TM.
\label{condition on X}
\end{align}
This is a stronger condition than the requirement of existence of  a timelike vector field in the Finslerian sense, namely,
\begin{align}
L(x,\mathcal{T}(x))=g_{(x,\mathcal{T}(x))}(\mathcal{T}(x),\mathcal{T}(x))<\,0.
\label{Finsler time orientation}
\end{align}
Given the vector field $\mathcal{X}$ such that \eqref{condition on X} holds good, for the fundamental tensor $g_{ij}(x,y)$ we define
\begin{align}
\tilde{g}_{(x,y)}(Y,Z):=\,g_{(x,y)}(Y,Z)-\,2\,\frac{g_{(x,y)}(\mathcal{X},Y)\,g_{(x,y)}(\mathcal{X},Z)}{g_{(x,y)}(\mathcal{X},\mathcal{X})}.
\label{lagrange space from g}
\end{align}
This is a generalization of the inverse transformation which associates to a given Riemannian structure, which always can be found for paracompact manifolds, a metric with Lorentzian signature \cite{Hawking Ellis 1973}.  Elementary algebraic manipulations leads to the expression
\begin{align*}
g_{(x,y)}(\mathcal{X},\mathcal{X})\,\tilde{g}_{(x,y)}(Y,Z)& =\,g_{(x,y)}(\mathcal{X},\mathcal{X})\,g_{(x,y)}(Y,Z)\\
& -\,2g_{(x,y)}(\mathcal{X},Y)\,g_{(x,y)}(\mathcal{X},Z).
\end{align*}
In the  case $\mathcal{X}=Y=Z$, it leads to the relation
\begin{align}
g_{(x,y)}(\mathcal{X},\mathcal{X})\,\tilde{g}_{(x,y)}(\mathcal{X},\mathcal{X})=\,-\left(g_{(x,y)}(\mathcal{X},\mathcal{X})\right)^2<\,0,\quad \forall \,(x,y)\in\,TM.
\label{relations timelike eigenvalue}
\end{align}
Therefore, the sign of $g_{(x,y)}(\mathcal{X},\mathcal{X})$ and the sign of $\tilde{g}_{(x,y)}(\mathcal{X},\mathcal{X})$ are opposite. On the other hand, for an orthonormal $Y(u)$ to $\mathcal{X}(u)$, then we have from the expression of $\tilde{g}$ that
\begin{align}
\tilde{g}_{(x,y)}(Y,Z):=\,g_{(x,y)}(Y,Z).
\label{spacelike component}
\end{align}
Hence if $Y,Z$ take values in an orthogonal basis for the bilinear form $g_{(x,y)}$, the same will be true for $\tilde{g}_{(x,y)}$.

The relations \eqref{relations timelike eigenvalue} and \eqref{spacelike component} imply that $\tilde{g}_{(x,y)}$ has positive signature, for each $u=(x,y)$.
 Note that, although $\tilde{g}$ is not the fundamental tensor of a positive definite Finsler space, but it is the fundamental tensor determining a Lagrange space \cite{MironHrimiucShimadaSabau:2002}.
\\
 {\bf 2}. It is remarkable that the averaging method can also be applied to positive definite Lagrange spaces, since the proposition \ref{propositionaveragemetric} and the formula \eqref{averagedmetric} applies directly to the fundamental $0$-homogenous tensor $g_{ij}(x,y)$, which is the fundamental object of a Lagrange space. We apply the averaging operation with the weight function $\psi=1$.
Then we have the following
\begin{proposicion}
 Let $(M,\tilde{g})$ be the Lagrange space given by the expression \eqref{lagrange space from g} and a vector field $\mathcal{X}$ such that the condition \eqref{condition on X} holds good. Then the components $\langle \,\tilde{g}_{ij}\rangle$ determine a Riemannian metric $\langle\tilde{g}\rangle$ on $M$.
 \label{Lagrange space}
\end{proposicion}

{\bf 3}. We can define now a metric with signature $n-1$ defined over $M$ from the Riemannian metric $\langle\tilde{g}\rangle$. We need to assume that there is a non-zero everywhere vector field $V\in\,\Gamma \,TM$, but we can assume that such a $V$ is indeed the vector field $\mathcal{X}$  introduced in adopting the condition \eqref{condition on X}. Therefore, we have the following
\begin{proposicion}
Let $(M,L)$ be a Finsler spaces with Lorentzian signature, $\mathcal{X}$ a vector field such that condition \eqref{condition on X} holds good. Then
the metric $\widetilde{\langle\tilde{g}\rangle}$ determined by the components whose expression is
\begin{align}
\widetilde{\langle\tilde{g}\rangle}(Y,Z):=\,\langle\tilde{g}\rangle_{(x,y)}(Y,Z)-\,2\,\frac{\langle\tilde{g}\rangle_{(x,y)}(\mathcal{X},Y)\,
\langle\tilde{g}\rangle_{(x,y)}(\mathcal{X},Z)}{\langle\tilde{g}\rangle_{x}(\mathcal{X},\mathcal{X})}
\label{Lorentzian space from g}
\end{align}
is a Lorentzian metric on $M$.
\label{Proposition Lorentzian space from q}
\end{proposicion}
This result suggests our definition of average Finsler spaces with Lorentzian signature,
\begin{definicion}
Given a Finsler spaces with Lorentzian signature $(M,L)$ and a vector field $\mathcal{X}$ as in \eqref{condition on X}, the average Finsler spaces with Lorentzian signature metric $\widetilde{\langle\tilde{g}\rangle}$ is given by the expression \eqref{Lorentzian space from g}.
\end{definicion}
It is direct that when the initial spacetime is Lorentzian ($g$ is defined over $M$), then $\mathcal{X}$ as in \eqref{condition on X} exists when the spacetime is time ordered \cite{GallegoPiccioneVitorio:2012}. Indeed, one has in the Lorentzian case that $\widetilde{\langle \tilde{g}\rangle}=\,\widetilde{\tilde{g}}=\,g$.

\subsection{Interpretation of the average Finsler spaces with Lorentzian signature structure}
We would like to remark that the condition \eqref{condition on X} is much stronger than the usual time orientation condition \eqref{Finsler time orientation}.
One can heuristically justify the condition \eqref{condition on X} if we think the Finsler spaces with Lorentzian signature structure $(M,L)$ as a small deformation from a time orientable Lorentzian structure. If $\mathcal{X}$ is the time orientation of the structure, then $L(x,\mathcal{X}(x))=g_{(x,\mathcal{X}(x))}(\mathcal{X}(x),\mathcal{X}(x))<\,0.$ Since $g(x,y)$ is homogeneous of degree zero in $y$-variable, it lives in sphere bundle $SM$, where each projective sphere is defined as the aggregate of equivalence classes
\begin{align*}
S_x M:\{[y]=\{y=\lambda\,y_0,\,\lambda\in\,\mathbb{R}^+,\,y_0\in \,N_x\}.
\end{align*}
Each sphere $S_xM$ is compact. Hence, if for each $x\in \,M$ the condition $L(x,\mathcal{X}(x))<\,0$ holds, then for small Finslerian perturbations of Lorentzian metrics the condition $g_{(x,y)}(\mathcal{X}(x), \mathcal{X}(x))<\,0$  with $y\in\,N$ is also to be expected.
 By accepting the condition \eqref{condition on X} we are accepting that $g_{(x,y)}$ and the Lorentzian metric $g_{(x,\mathcal{X}(x))}$ do not differ too much.

 A fundamental difference of the average Finsler structures with Lorentzian signature structure respect to the analogous operation for Finsler metrics with positive signature is that in the Lorentzian case, $\tilde{g}$ depends on the timelike vector field $\mathcal{X}\in\,\Gamma TM$ and on the timelike vector field $V$ defining  the metric $\widetilde{\langle\tilde{g}\rangle}$, although in our construction we have taken $V=\mathcal{X}\in\,\Gamma \,TM$. This makes the metric $\widetilde{\langle\tilde{g}\rangle}$ to depend on the vector field $\mathcal{X}$.

\section{On the average connection}
Let $(M,F)$ be a positive definite Finsler space. As it has been discussed in \cite{Ricardo05}, there is defined the Chern connection $\nabla$ and its average connection $\langle \nabla\rangle$ as defined by {\it Theorem} \ref{averagedconnection}. It can be shown that the average connection $\langle\nabla \rangle$ is an affine connection on $M$ \cite{Ricardo05}.
The construction of the average Chern connection can be extended to the case of Finsler spaces with Lorentzian signature $n-1$. Applying the same procedure as in the case of positive definite Finsler space, after implementing the adequate changes in the measure used in the averaging method when passing from positive definite to signature $n-1$, the connection coefficients in a particular holonomic basis for $T_xM$ for the average metric are given by the expression
\begin{align}
\langle \Gamma^i_{jk}\rangle_{\psi} (x):=\frac{1}{vol(I_x)} \int _{\Sigma_x }\,|\psi|^2(x,y)\,\Gamma^i_{jk}(x,y)\,dvol_x,
\label{average connection 1}
\end{align}
where the connection coefficients $\Gamma^i_{jk}(x,y)$ are the connection coefficients of the Chern connection.

 The submanifold $\Sigma_x\subset \,N_x$ is non-compact. Therefore, in order  to ensure the convergence of the integral, the measure $|\psi|^2(x,y)$ must converge in the limit $y$-coordinates going to infinity conveniently. In general, it is a difficult task to find an universal measure $|\psi|^2$ with such characteristics. Hence we need to assume that $\psi$ depends on the particular details of the given function $L$.
 
 Another possible rout to avoid the problems on the non-compactness of the integration domain $\Sigma_x$  is to define the averaging operation of metrics as an integral operation on a subset of the projective sphere bundle $S_xM$. Since $S_xM$ is compact, there is no need of an special convergence property for the measure. Therefore, we can take $|\psi(x,y)|=1$.
 On each projective sphere bundle $S_x M$ there is defined an induced volume form $dvol'_x$ determined by the volume form $dvol_x$ defined on $N_x$. Let us consider the canonical projection $p:N_x\to  S_xM$ and let us note the relation  $\dim (N_x)=\dim (S_xM)-1$. Then $dvol'_x$ is defined as the unique volume form on $S_xM$ such that the pull-back $p^*dvol'_x$ is equal to the contraction $\,\l_x\lrcorner \,dvol_x$, where $\l_x=y^k\frac{\partial}{\partial y^k}|_x$.
Then one can define the following average
\begin{align}
\langle \Gamma^i_{jk} \rangle(x):=\frac{1}{vol(I_x)} \int _{S_x M}\,\Gamma^i_{jk}(x,y)\,dvol'_x.
\label{average connection 2}
\end{align}
\begin{ejemplo}
Of particular interest are Finsler spaces with Lorentzian signature where the connection coefficients of the Chern connection do not depend upon the $y$-variable. Examples of such spacetimes are Berwald type spacetimes \cite{Ricardo 2017b}. Fixed a normal coordinate system on $TM$, let us consider the following sub-manifold of $T_xM$,
\begin{align*}
P_x(r)=\{y\in\,T_xM\,s.t.\, |y^1|\leq r\}.
 \end{align*}
  The integral operation is the regularized integral \eqref{average connection 1},
\begin{align*}
\langle \Gamma^i_{jk}\rangle_{{\psi}_r}(x):=\frac{1}{vol(P_x(r))} \int _{P_x(r) }\,|\psi|^2(x,y)\,\Gamma^i_{jk}(x)\,dvol_x=\,\Gamma^i_{jk}(x),
\end{align*}
where
\begin{align*}
{vol(P_x(r))}=\, \int _{P_x(r) }\,|\psi|^2(x,y)\,dvol_x.
\end{align*}
In the limit $r\to \infty$ the relation
\begin{align*}
\langle \Gamma^i_{jk}\rangle_{{\psi}_r}(x)\to\,\langle \Gamma^i_{jk}\rangle_{\psi} (x), \quad i,j,k=1,...,n
\end{align*}
holds good, showing that the limit $r\to\,+ \infty$ in the integral operation is well defined and equal to $\Gamma^i_{jk}(x)$.
\end{ejemplo}
It is easy to see that, as in the positive definite case, the torsion of the averaged Chern connection given by either the connection coefficients \eqref{average connection 1} or by \eqref{average connection 2}, is zero, However, in general the averaged Chern connection is not the Levi-Civita connection of the average metric $\widetilde{\langle \tilde{g}\rangle}$, even when the Finsler spaces with Lorentzian signature is of Berwald type.

\section{Conclusion and outlook}
We have shown how to extend the averaging operations of metrics and connections from positive definite Finsler structures to certain Finsler structures with signature $n-1$. In the case of the metric structure, the main requirement added for the construction to work is the existence of a vector field $\mathcal{X}$ such that the condition \eqref{condition on X} is fulfilled. This condition has heuristically motivated. However, the problem of defining a natural average operation in the general case, remains still an open problem.

We have also discussed the problem of averaging the Chern connection in the setting of metrics with signature $n-1$. In this case, the procedure does not need the existence or not of the vector field $\mathcal{X}$ constrained by \eqref{condition on X}. This problem admits a general solution following the methodology described in this paper, except for some issues of convergence due to the non-compactness of the integration domain $\Sigma_x$.

Since the procedures for averaging the metric and the connection are rather different, one should not expect a direct translation of the fundamental characterization of Berwald spaces as given in \cite{RF}, namely, that a space is of Berwald type iff the averaged Chern connection is the Levi-Civita connection of the average metric.
It is an interesting problem, due to the relevance of Berwald spacetimes for physics, to know if there is a definition of the average connection in for signature $n-1$  such that the characterization given in \cite{RF}  still holds.

\end{document}